\documentclass[10pt, a4paper]{amsart}
\usepackage{amsfonts}
\usepackage{amssymb}
\usepackage{amsthm}
\usepackage{amsmath}
\usepackage{enumitem}
\usepackage{pifont}
\usepackage{float}
\usepackage{pgfplots}
\usepackage{xfrac}
\usepackage{anyfontsize}
\usepackage{color}
\usepackage{changepage}
\usepackage{url}
\usepackage{bm}
\usepackage{datetime}
\usepackage{faktor}
\usepackage{tikz}
\usepackage[abs]{overpic}
\usepackage{tikz-cd}
\usepackage{forest}
\usetikzlibrary{calc, shapes, decorations.pathmorphing}
\usepackage{graphicx}
\usepackage{csquotes}

\usepackage{hyperref}
\usepackage{mathtools}
\allowdisplaybreaks

\pgfdeclarelayer{nodelayer}
\pgfdeclarelayer{edgelayer}
\pgfsetlayers{nodelayer,edgelayer,main}

\def \proof {\noindent\textit{Proof.}\hspace{5mm}}
\def \bs {$\hfill \blacksquare$\\}

\DeclareMathOperator{\red}{red}

\DeclareMathOperator{\charac}{char}

\DeclareMathOperator{\Z}{\mathbb{Z}}

\DeclareMathOperator{\Q}{\mathbb{Q}}

\DeclareMathOperator{\p}{\mathfrak{p}}
\DeclareMathOperator{\q}{\mathfrak{q}}

\definecolor{supergreen}{RGB}{0, 170, 0}
\definecolor{superred}{RGB}{170, 0, 0}
\definecolor{darkred}{RGB}{100, 0, 0}
\definecolor{customcolour}{RGB}{0, 128, 128}

\newcommand{\ThE}{\mbox{\ttfamily Th}_{\exists}}

\newcommand{\fs}{\mathbb{F}_q [t]}
\newcommand{\fr}{\mathbb{F}_q (t)}
\newcommand{\iso}{\cong}
\newcommand{\fin}{\mathbb{F}_q \! \left(\!\left({}^1{\mskip -5mu/\mskip -3mu}_t\right)\!\right)}
\newcommand{\fract}{{}^1{\mskip -5mu/\mskip -3mu}_t}

\newcommand{\defeq}{\mathrel{\mathop:}=}

\newcommand{\st}{\mbox{ : }}

\newtheorem*{thm*}{Theorem}
\newtheorem{thm}{Theorem}[section]

\newtheorem{cor}[thm]{Corollary}
\newtheorem{lem}[thm]{Lemma}

\newtheoremstyle{case}{}{}{}{}{}{:}{ }{}
\theoremstyle{case}

\theoremstyle{definition}
\newtheorem{definition}[thm]{Definition}

\pgfplotsset{compat=1.14}

\begin{document}
 
\title{A New Universal Definition of $\fs$ in $\fr$}
\author{Brian Tyrrell}
\thanks{\textit{2010 Mathematics Subject Classification: } 12L12 (primary) and 03C60 (secondary).}
\address{Mathematical Institute, Woodstock Road, Oxford OX2 6GG.}
\email{brian.tyrrell@maths.ox.ac.uk}

\begin{abstract}
This paper gives a universal definition of $\fs$ in $\fr$ using 89 quantifiers, more direct than those that exist in the current literature. The language $\mathcal{L}_{\mbox{\tiny rings}, t}$ we consider here is the language of rings $\{0, 1, +, -, \cdot\}$ with an additional constant symbol $t$. We then modify this definition marginally to universally define $\fs$ in $\fr$ \textit{without} parameters, using 90 quantifiers. We assume throughout that $\charac(\mathbb{F}_q) \neq 2$.
\end{abstract}

\maketitle

\section{Introduction}
One motivation for definability questions (such as determining a universal definition of $\fs$ in $\fr$, or an existential definition of $\Z$ in $\Q$) stems from David Hilbert's famous list of 23 problems, published in 1900 \cite{hilbert}. In particular, his tenth problem (H10) requests to prove $\ThE(\Z)$ is decidable. We know now after the work of Davis, Putnam, Robinson, and Matiyasevich \cite{dpr, matiya2} that this theory is in fact undecidable. However, as is often the case in mathematics, the answer to a good question itself raises more questions than answers available. The `natural' generalisation of H10 is to determine the decidability of $\ThE(\Q)$, and this question still remains open. There are two paths, amongst others, before us in the quest to answer ``H10 over $\Q$'' -- one path's approach is via the definability of $\Z$ in $\Q$, the other is via H10 over other rings and fields. The former approach is useful as follows: if one had an existential definition of $\Z$ in $\Q$, then the undecidability of $\ThE(\Q)$ follows from the undecidability of $\ThE(\Z)$. The latter approach is more philosophical; if one understood the behaviour of H10 over different rings and fields, one could possibly gain a deeper insight into the problem and use this to solve H10 over $\Q$, and go further with further generalisations such as H10 over $K$ or $\mathcal{O}_K$ where $K$ is a number field.

This paper lies in the middle between these two paths. A significant addition to the definability approach of solving H10 over $\Q$ was made in 2016 by Koenigsmann \cite{defzinq}, who provided a universal definition (later improved by Daans \cite{daans,daansnew}) and a $\forall\exists$-definition of $\Z$ in $\Q$, the latter only using \emph{one} universal quantifier. For the latter more philosophical path it is worth noting that both $\ThE(\fs)$ and $\ThE(\fr)$ are undecidable (in $\mathcal{L}_{\mbox{\tiny rings}, t}$) \cite{denef, pheidas, videla} which seems almost to suggest (by the \emph{function field analogy}) that a more complete understanding of H10 made in the function field context would be useful for determining H10 over $\Q$. This paper is not the first to explore definability questions in function fields; Eisentr\"{a}ger \&\ Morrison \cite{eisen} adapted Park's \cite{park} universal definition of $\mathcal{O}_K$ in $K$ from number fields to function fields, and this definition was vastly simplified by Daans \cite{daansnew} who in fact provided a universal definition of the ring of $S$-integers $\mathcal{O}_S$ in a global field $G$ where $S$ is a finite nonempty set of primes of $G$. In \cite{daans} there is also a shorter, more easily attained universal definition of $\Z$ in $\Q$ than Koenigsmann (\cite[Theorem 4.3.3]{daans}), and it is from this theorem that the main result of the paper sprang.

\begin{thm*}
Assume $\charac(\mathbb{F}_q) \neq 2$. There is a universal definition of $\fs$ in $\fr$ given by 89 quantifiers, and a universal $\emptyset$-definition given by 90 quantifiers.
\end{thm*}

The essence of \cite[Theorem 4.3.3]{daans} can be summarised as follows: the main goal of the theorem is to find an existentially defined set of conditions $E$ on parameters $a,b$ such that
\begin{enumerate}
    \item If $a,b$ satisfies $E$ this forces $\Delta_{a,b}^{\circ} \defeq \Delta_{a,b} \cap (\mathbb{P}(a) \cup \mathbb{P}(b)) = \Delta_{a,b} \setminus \{\infty\}$, where
    \begin{gather*}
        \Delta_{a,b} \defeq \{ p \in \mathbb{P}\cup\{\infty\} \st H_{a,b} \otimes \Q_p \not\iso M_2(\Q_p) \mbox{ (i.e.\ does not split)}\},\\
        \mathbb{P}(a) \defeq \{p \in \mathbb{P} \mbox{ : } v_p(a) \mbox{ is odd}\},
    \end{gather*}
    and $H_{a,b}$ is the quaternion algebra $\Q \cdot 1 \oplus \Q \cdot \alpha \oplus \Q \cdot \beta \oplus \Q \cdot \alpha \beta$ with $\alpha^2 = a, \beta^2 = b, \alpha \beta = -\beta \alpha$.
    \item If $a,b$ satisfy $E$, then $(a, b)_{\infty} = -1$, where $(a,b)_{p}$ is the \emph{Hilbert symbol}
    $$(a, b)_p = 1 \mbox{ when } p \not\in \Delta_{a,b}, \quad (a, b)_p = -1 \mbox{ when } p \in \Delta_{a,b}.$$    
    This will ensure, by Hilbert Reciprocity, $\Delta_{a,b}^{\circ}$ is always nonempty.
    \item For each prime $p$, one can find $a, b$ satisfying $E$ such that $\Delta_{a,b}^{\circ} = \{p\}$. Equivalently, there exist $a, b$ satisfying $E$ such that
    $$(a,b)_{p} = -1 \mbox{ and } (a,b)_{q} = 1 \mbox{ for all primes }q \neq p.$$
\end{enumerate}

Then we obtain a universal definition
\begin{itemize}
    \item[(4)] $t \in \Z \quad \Leftrightarrow \quad \forall a, b \in \Q \mbox{ }((a,b) \not \in E \lor t \in \widetilde{R_{a,b}})$,
\end{itemize}
where $\widetilde{R_{a,b}}$ is a universally defined union of localisations of $\Z$ at primes $p \in \Delta_{a,b}^{\circ}$.\\

To accommodate the fact that all primes of $\fr$ are, in some sense, ``finite'' (nonarchimedian) we will have to modify (1) in order for $\widetilde{R_{a,b}}$ to have a universal definition in this setting. We shall find a new set of existentially defined conditions $D$ on parameters $a,b$ such that
\begin{itemize}
    \item[(1')] If $a,b$ satisfies $D$ this forces $\Delta_{a,b}^{\circ} = \Delta_{a,b}$.
    \item[(2)] If $a,b$ satisfy $D$, then $(a, b)_{\infty} = -1$, where $\infty$ is the prime of $\fr$ corresponding to the valuation $-\deg$, and $\fr_{\infty} = \fin$.
    \item[(3)] For each prime $\mathfrak{p}$, one can find $a, b$ satisfying $D$ such that $\Delta_{a,b}^{\circ} = \{\mathfrak{p, \infty}\}$.
\end{itemize}
Then we will obtain a universal definition as follows: writing $(\fr)_{\mathfrak{p}}$ for $\fr$ localised at a prime $\p$,
\begin{itemize}
    \item[(4)] $t \in \fs \cup \left(\fs\right)_{\infty} \quad \Leftrightarrow \quad \forall a, b \in \fr \mbox{ }((a,b) \not \in D \lor t \in \widetilde{R_{a,b}})$,
\end{itemize}
from which a universal definition of $\fs$ in $\fr$ can be quickly obtained.

At the time of writing this led to the shortest (in number of quantifiers) known universal definition of $\fs$ in $\fr$, however using some intricate quaternion algebra theory and deep class field theory, Daans \cite{daansnew} proves there is a universal definition of $\fs$ in $\fr$ requiring only a breathtaking 65 quantifiers.

Let us begin our definition by first determining $D$.

\bigskip
\section{A New Universal Definition}\label{shorteruni}

We will assume that $\charac(\mathbb{F}_q) \neq 2$ (necessary for \emph{Lemma \ref{second} \&\ Theorem \ref{omgproof}}). We first need the following characterisation of nonsquares of $\fin$:



\begin{lem}\label{first}
Any nonsquare of $\fin$ is of the form $\fract c^2, fc^2$, or ${}^f{\mskip -5mu/\mskip -3mu}_t c^2$ where $c \in \fin$ and $f \in \mathbb{F}_q$ is a nonsquare.
\end{lem}

\proof
This is an application of Hensel's Lemma exactly. 
\bs

\begin{lem}\label{second}
The quaternion algebra 
$$H_{f, {}^g{\mskip -5mu/\mskip -3mu}_t }(\fin) = \fin \cdot 1 \oplus \fin \cdot \alpha \oplus \fin \cdot \beta \oplus \fin \cdot \alpha \beta$$
with multiplication defined by $\alpha^2 = f, \beta^2 = {}^g{\mskip -5mu/\mskip -3mu}_t, \alpha\beta = -\beta\alpha$, is nonsplit, where $f,g \in \mathbb{F}_q^{\times}$ and $f$ is a nonsquare.
\end{lem}

\proof
Using \cite[XIV.3.8]{serre79}, for a $\p$-adic unit $a$,
$$(a, b)_{\p} = -1 \quad\Leftrightarrow\quad v_{\p}(b) \mbox{ is odd and } \red_{\p}(a) \mbox{ is a nonsquare of } \mathbb{F}_{\p}.$$

Thus $(f, {}^g{\mskip -5mu/\mskip -3mu}_t )_{\infty} = -1$ if and only if $v_{\infty}({}^g{\mskip -5mu/\mskip -3mu}_t ) = 1$ is odd and $f \in \mathbb{F}_q^{\times}$ is a nonsquare (as it was chosen to be). Hence $H_{f, {}^g{\mskip -5mu/\mskip -3mu}_t}(\fin)$ is nonsplit. Note this also means $H_{{}^f{\mskip -5mu/\mskip -3mu}_t, g}(\fin)$ is nonsplit too.
\bs

We adopt the following piece of notation: if $\p = (f(t))$ is a prime of $\fs$ (where $f(t)$ is a monic and irreducible polynomial) then the residue field of $\fr_{\p} = \fr_{f(t)}$ under the $\p$-adic valuation is denoted $\mathbb{F}_{f(t)}$ and is isomorphic to the set of polynomials of $\fs$ of degree strictly less than $\deg(f)$. The residue map $\fs_{f(t)} \rightarrow \mathbb{F}_{f(t)}$ is denoted $\red_{f(t)}$. We will make use of the Legendre symbol, which in this context is defined as:

\begin{definition}
Let $f(t) \in \fs$ be a prime (that is, the monic and irreducible polynomial corresponding to the principal prime ideal $\p$) and $g(t) \in \fs$, where $f(t) \nmid g(t)$. Then
\begin{align*}
    \left(\tfrac{g(t)}{f(t)}\right) \defeq
    \begin{cases}
    1 &\mbox{if }\red_{f(t)}(g(t)) \mbox{ is a square of } \mathbb{F}_{f(t)},\\
    -1 &\mbox{if }\red_{f(t)}(g(t)) \mbox{ is a nonsquare of } \mathbb{F}_{f(t)}.
    \end{cases}
\end{align*}
\end{definition}
\vspace{1mm}
\begin{lem}\label{gotit}
Let $f(t) \in \fs$ be a prime and $g \in \mathbb{F}_q$ be nonsquare. If $\deg(f)$ is odd, then $\left(\tfrac{g}{f(t)}\right) = -1$. If $\deg(f)$ is even, then $\left(\tfrac{g}{f(t)}\right) = 1$.
\end{lem}

\proof
This follows from \cite[Prop.\ 3.2]{rosen}.\bs

\vspace{-1mm}
\begin{lem}\label{finally}
Given a prime $f(t) \in \fs$ and $g \in \mathbb{F}_q$ nonsquare, one can choose $d(t)$ a prime of $\fs$ of opposite parity in degree to $f(t)$ such that $\red_{f(t)}(g d(t))$ is a nonsquare of $\mathbb{F}_{f(t)}$. Moreover, this choice can be made independent of $g$.
\end{lem}

\proof
By Dirichlet's Theorem on primes in arithmetic progressions there are infinitely many primes equivalent to $c(t)$ mod $f(t)$ for any $c(t) \in \mathbb{F}_{f(t)}$. Moreover, for $N$ large enough, there is a prime of degree $N$ in this arithmetic progression \cite[Theorem 4.8]{rosen}. 

Therefore if $f(t)$ is of odd degree then we can choose $d(t)$ to be monic, irreducible, of even degree and $d(t) \equiv c(t)^2 \mod f(t)$, where $c(t) \not\equiv 0 \mod f(t)$. If $f(t)$ has even degree then we can choose $d(t)$ to be monic, irreducible, of odd degree and $d(t) \equiv c(t) \mod f(t)$ where $c(t) \in \mathbb{F}_{f(t)}$ is a nonsquare. Then $\red_{f(t)} (gd(t))$ is a nonsquare of $\mathbb{F}_{f(t)}$, according to \emph{Lemma \ref{gotit}}.
\bs

These lemmata will contribute to the next result. First, we introduce more notation. For $c \in \fr$, let $\phi(c)$ denote ``$c$ as an element of $\fin$ is a square''. An equivalent statement, by \emph{Lemma \ref{first}}, is that the degree of $c$ is even and \big(writing $c = \tfrac{a_m t^m + \dots + a_0}{b_n t^n + \dots + b_0}$\big) the \emph{leading coefficient} $\tfrac{a_m}{b_n}$ is a square. Let $z \in \mathbb{F}_q$ be a nonsquare and let $\psi(a,b)$ denote
\begin{align*}
    \exists c, d \Big(&\mbox{``$c$ and $d$ are of opposite parity in degree''}\\
    &\land \Big[ \big\{ \phi(c) \land a = z c \land b^{q-1} = d^{q-1} \big\} \lor \big\{ \phi(d) \land b = z d \land a^{q-1} = c^{q-1} \big\} \Big]\Big).
\end{align*}

Finally define
\begin{definition}\label{D}
$D \defeq \{(a, b) \in \fr \times \fr \mbox{ : } \psi(a,b)\}.$
\end{definition}

The complicated choice of $\psi(a,b)$ will be justified in the upcoming theorem.

\begin{thm}\label{omgproof}
We have
\begin{equation*}\label{omg}
\fs \cup \left(\fs\right)_{\infty} = \bigcap_{(a,b) \in D} \widetilde{R_{a,b}},  
\end{equation*}
where
$$R_{a,b} \defeq \bigcap_{\p \in \Delta_{a,b} \cap (\mathbb{P}(a) \cup \mathbb{P}(b))} (\fs)_{\p}, \qquad \widetilde{R_{a,b}} = \bigcup_{\p \in \Delta_{a,b} \cap (\mathbb{P}(a) \cup \mathbb{P}(b))} (\fs)_{\p}.$$

\end{thm}

\proof
To begin, consider the set of primes $\Delta_{a,b}$ in more detail. 
\begin{align*}
    \p \in \Delta_{a,b} \quad&\Leftrightarrow\quad (a,b)_{\p} = -1 \\
    &\Leftrightarrow\quad \bigg((-1)^{v_{\p}(a) v_{\p}(b)} \red_{\p} \left(\frac{a^{v_{\p}(b)}}{b^{v_{\p}(a)}}\right) \bigg)^{\tfrac{\# \mathbb{F}_{\p} - 1}{2}} = -1, \mbox{ by \cite[XIV.3.8]{serre79}}.
\end{align*}

Assume for the purpose of contradiction that $\p \not\in \mathbb{P}(a) \cup \mathbb{P}(b)$: then $v_{\p}(a)$ and $v_{\p}(b)$ are both even. Assume one of them is nonzero.\footnote{If $v_{\p}(a) = v_{\p}(b) = 0$, then $c = 1$ and $\red_{\p}(c)^{\tfrac{\#\mathbb{F}_q - 1}{2}} = 1$; a contradiction too.}
\begin{align*}
&\Leftrightarrow\quad \bigg(\red_{\p}(c)^2\bigg)^{\tfrac{\# \mathbb{F}_{\p} - 1}{2}} = -1, \mbox{ where }c=\tfrac{a^{v_{\p}(b)/2}}{b^{v_{\p}(a)/2}}\\
    &\Leftrightarrow\quad \red_{\p}(c)^{\#\mathbb{F}_{\p} - 1} = -1,
\end{align*}
however $\red_{\p}(c)$ must satisfy the equation $x^{\#\mathbb{F}_{\p}} = x$ of a finite field; with our assumption of a noneven characteristic, we have a contradiction. Thus
$$\Delta_{a,b}^{\circ} \quad= \quad \Delta_{a,b} \cap (\mathbb{P}(a) \cup \mathbb{P}(b)) \quad = \quad \Delta_{a,b}.$$

We will now prove $\Delta_{a,b}$ is nonempty for $(a,b) \in D$: any nonsquare of $\fr_{\infty} = \fin$ is of the form $\fract c^2, fc^2$, or ${}^f{\mskip -5mu/\mskip -3mu}_t c^2$ for $c \in \fin$ and $f \in \mathbb{F}_q$ a nonsquare, by \emph{Lemma \ref{first}}. For $(a, b) \in D$ considered as elements of $\fin$, there are at most 9 possible classes for $(a,b)$ modulo squares of $\fin$:

\begin{center}
\begin{tabular}{c c c}
$(\fract, \fract)$ & $(\fract, f)$ & $(\fract, {}^f{\mskip -5mu/\mskip -3mu}_t )$   \\
$(f, \fract)$ & $(f, g)$ & $(f, {}^g{\mskip -5mu/\mskip -3mu}_t )$   \\
$({}^f{\mskip -5mu/\mskip -3mu}_t , \fract)$ & $({}^f{\mskip -5mu/\mskip -3mu}_t , g)$ & $({}^f{\mskip -5mu/\mskip -3mu}_t , {}^g{\mskip -5mu/\mskip -3mu}_t )$ 
\end{tabular}
\end{center}
for $f, g \in \mathbb{F}_q$ nonsquares. However out of these possible scenarios, only four are allowed by choice of $a$ and $b$: $(z, {}^g{\mskip -5mu/\mskip -3mu}_t ), ({}^g{\mskip -5mu/\mskip -3mu}_t, z), ({}^1{\mskip -5mu/\mskip -3mu}_t, z)$ and $(z, {}^1{\mskip -5mu/\mskip -3mu}_t)$. By the rules of quaternionic bases we conclude $H_{a,b}(\fin)$ is nonsplit for any such $a, b$ if $H_{z, {}^g{\mskip -5mu/\mskip -3mu}_t }(\fin)$ and $H_{z, {}^1{\mskip -5mu/\mskip -3mu}_t }(\fin)$ are nonsplit. However by \emph{Lemma \ref{second}} we know these are nonsplit.

This demonstrates that if $(a,b) \in D$, then $\infty \in \Delta_{a,b}$. As well as this, by Hilbert Reciprocity we conclude the quaternion algebra given by $(a,b)$ must be nonsplit at some non-infinite prime too, meaning $\Delta_{a,b} \setminus \{\infty\}$ is nonempty. This allows us to conclude $\fs \cup \left(\fs\right)_{\infty} \subseteq \widetilde{R_{a,b}}$ for each $(a,b) \in D$, therefore 
$$\fs \cup \left(\fs\right)_{\infty} \subseteq \bigcap_{(a,b) \in D} \widetilde{R_{a,b}}.$$

We will now prove the reverse inclusion. Consider the prime ideals of $\fs$; these are principal ideals $\p = (f(t))$ with $f(t) \in \fs$ a monic and irreducible polynomial.

Set $a = z f(t)$ and $b = z d(t)$ according to \emph{Lemma \ref{finally}}. By this choice of $a$ and $b$, $(a,b)_{\p} = -1$ as $v_{\p}(a)$ is odd and $\red_{\p}(b)$ is a nonsquare of $\mathbb{F}_{\p}$. Also, if $\q$ is a prime such that $\q \neq \p$, $\q \neq \infty$, it follows that $v_{\q}(a) = 0$ and $b$ is either a $\q$-unit (in which case $(a,b)_{\q} = 1$) or $\q = (d(t))$. In this case,
\begin{align*}
    (a,b)_{d(t)} &= \left( (-1)^{v_{d(t)}(a) v_{d(t)}(b)} \red_{d(t)} \left(\tfrac{a^{v_{d(t)}(b)}}{b^{v_{d(t)}(a)}}\right) \right)^{\tfrac{q^{\deg d} - 1}{2}}\\
    &= \red_{d(t)}(z f(t))^{\tfrac{q^{\deg d} - 1}{2}}\\
    &= \left(\tfrac{zf(t)}{d(t)}\right) = \left(\tfrac{z}{d(t)}\right) \left(\tfrac{f(t)}{d(t)}\right).
\end{align*}
By the law of Quadratic Reciprocity (cf.\ \cite[Theorem 3.3]{rosen}),
\begin{equation*}
    \left(\tfrac{d(t)}{f(t)}\right)\left(\tfrac{f(t)}{d(t)}\right) = (-1)^{\tfrac{q-1}{2} \deg f \deg d} = 1,
\end{equation*}
as $f$ and $d$ have opposite parity in degree (and $q$ is not a power of 2). Consider the following two cases.
\begin{itemize}
    \item[\textbf{Case 1:}] $f$ has odd degree. Then $\left(\tfrac{d(t)}{f(t)}\right) = 1$ by \emph{Lemma \ref{finally}}, meaning $\left(\tfrac{f(t)}{d(t)}\right) = 1$. Also $\left(\tfrac{z}{d(t)}\right) = 1$ by \emph{Lemma \ref{gotit}}, so
    $$(a,b)_{d(t)} = \left(\tfrac{z}{d(t)}\right) \left(\tfrac{f(t)}{d(t)}\right) = (1)(1) = 1.$$
    
    \item[\textbf{Case 2:}] $f$ has even degree. Then $\left(\tfrac{d(t)}{f(t)}\right) = -1$ by \emph{Lemma \ref{finally}}, meaning $\left(\tfrac{f(t)}{d(t)}\right) = -1$. Also $\left(\tfrac{z}{d(t)}\right) = -1$ by \emph{Lemma \ref{gotit}}, so
    $$(a,b)_{d(t)} = \left(\tfrac{z}{d(t)}\right) \left(\tfrac{f(t)}{d(t)}\right) = (-1)(-1) = 1.$$
\end{itemize}

In either case, we conclude $(a,b)_{d(t)} = 1$. So by choice of $a$ and $b$, $\p$ and naturally $\infty$ are the only primes at which the algebra $H_{a,b}(\fr_{\p})$ is nonsplit. Moreover by design $(a,b) \in D$, so $\Delta_{a,b} = \{\p, \infty\}$ and
$$\bigcap_{(a,b) \in D} \widetilde{R_{a,b}} \subseteq \left( \bigcap_{\p \neq \infty} (\fs)_{\p}\right) \cup (\fs)_{\infty} = \fs \cup (\fs)_{\infty},$$
as required.
\bs

We will show now that $D$ of \emph{Definition \ref{D}} is existentially definable.

\begin{lem}\label{therem}
Let $z \in \mathbb{F}_q$ be a nonsquare and let $\psi(a,b)$ denote
\begin{align*}
    \exists c, d \Big(&\mbox{``$c$ and $d$ are of opposite parity in degree''}\\
    &\land \Big[ \big\{ \phi(c) \land a = z c \land b^{q-1} = d^{q-1} \big\} \lor \big\{ \phi(d) \land b = z d \land a^{q-1} = c^{q-1} \big\} \Big]\Big).
\end{align*}
Then $\psi(a,b)$ is equivalent to an existential formula.
\end{lem}

\proof
For $c \in \fr$, consider $\phi(c)$: ``$c$ as an element of $\fin$ is a square''. Quantifying over $\fr$, this is captured by
\begin{align*}
    &\exists f \Big( \deg(c) = \deg(f^2) \land \exists g \big( \deg(c) > \deg(g) \land c = f^2 + g \big) \Big) \\
    &\Leftrightarrow\mbox{ }\exists f \big( \deg(c) = \deg(f^2) \land \deg(c) \geq \deg(t(c - f^2)) \big)    \\
    &\Leftrightarrow \mbox{ } \exists f \big( \deg(f^2) \geq \deg(t(c - f^2)) \big)
\end{align*}
What if we additionally wanted to say ``and $d$ is of odd degree''? This would be
\begin{align*}
    &\mbox{\hspace{1.15mm}}\exists f \big( \deg(f^2) \geq \deg(t(c - f^2)) \big) \land \exists h \big( \deg(f^2) = \deg(th^2 d) \big) \\
    &\Leftrightarrow \exists f, h \big( \deg(t h^2 d) \geq \deg(t(c - f^2)) \big) \land \exists g \big( \deg(f^2) > \deg(g) \land f^2 = k \cdot t h^2 d + g \big)\\
    &\mbox{    for some } k \in \mathbb{F}_q, \\
    &\Leftrightarrow \exists f, h \big( \deg(t h^2 d) \geq \deg(t (c - f^2)) \land \deg(f^2) \geq \deg(t(f^2 - k\cdot th^2 d)) \big)\\
    &\mbox{    for some } k \in \mathbb{F}_q.
\end{align*}
Let $\chi(c,d)$ denote
\begin{align*}
    \exists f, h \big( \deg(t h^2 d) \geq \deg(t (c - f^2)) \land \big\{ &\deg(f^2) \geq \deg(t(f^2 - th^2 d))\\
    &\lor \deg(f^2) \geq \deg(t(f^2 - z\cdot th^2 d))\big\} \big).
\end{align*}

Then, by the above argument and \emph{Lemma \ref{first}}, ``the degree of $c$ is even, the degree of $d$ is odd, and the leading coefficient of $c$ is a square'' is represented by this formula. Therefore $\psi(a,b)$ is equivalent to
\begin{equation}\label{thetag}
\chi(za, b) \lor \chi(z b, a).    
\end{equation}

The formula ``$\deg(A) \geq \deg(B)$'' is equivalent to ``$v_{\infty}(\tfrac{B}{A}) \geq 0$''. By \cite[Theorem 3.1]{eisenint}, the set $\{z \in K \mbox{ : } v_{\p}(z) \geq 0\}$ is existentially definable (and requires 9 quantifiers to define), therefore $\psi(a,b)$ is indeed equivalent to an existential formula and moreover requires $2 + 9 + 9 = 20$ quantifiers according to (\ref{thetag}).\bs

\begin{cor}\label{thecor}
There is a universal definition of $\fs$ in $\fr$ given by 89 quantifiers.
\end{cor}

\proof
By \emph{Theorem \ref{omgproof}}, we have
\begin{equation}\label{uni}
f(t) \in \fs \cup \left(\fs\right)_{\infty}\quad \Leftrightarrow \quad\forall a, b\left( (a,b) \not\in D \lor f(t) \in \widetilde{R_{a,b}}\right).    
\end{equation}

By \cite[Lemma 3.19]{eisen}, $\widetilde{R_{a,b}}$ is universally defined, hence as $D$ is existentially defined, (\ref{uni}) is indeed a universal formula for $\fs \cup \left(\fs\right)_{\infty}$. Denote this formula by $\Phi$. Recall that the number of quantifiers needed to define $\widetilde{R_{a,b}}$ is one more than is required to define its Jacobson radical. By \cite[Lemma 3.17]{park}, $J(\widetilde{R_{a,b}})$ has an explicit description of
\begin{align*}
    \{0\} \cup \{ x \in \fr^{\times} \mbox{ : } \exists y_1, y_2 \mbox{ s.t. } &y_1, x - y_1 \in a \cdot \fr^2 \cdot T_{a,b}^{\times} \cap (1 - \fr^2 \cdot T_{a,b}^{\times})\\
    &y_2, x - y_2 \in b \cdot \fr^2 \cdot T_{a,b}^{\times} \cap (1 - \fr^2 \cdot T_{a,b}^{\times})\},
\end{align*}
where $T_{a,b} = S_{a,b} + S_{a,b}$ and
$$S_{a,b} = \{2 x_1 \in \fr \mbox{ : } \exists x_2, x_3, x_4 \mbox{ s.t. } x_1^2 - a x_2^2 - b x_3^2 + ab x_4^2 = 1 \}.$$
The set $S_{a,b}$ requires three quantifiers to define, hence $T_{a,b}$ is defined by 7 quantifiers, as is $T^{\times}_{a,b}$ (by \cite[Lemma 4.3]{daansnew}). Thus $J(\widetilde{R_{a,b}})$ requires at most 66 quantifiers to define it. Thus the number of universal quantifiers need to define $\fs \cup \left(\fs\right)_{\infty}$ in $\fr$ using (\ref{uni}) is at most $2 + 20 + 66 + 1 = 89$. 

What about the definition of $\fs$? This is simply
$$f(t) \in \fs \quad \Leftrightarrow \quad \Phi(f(t)) \land \left( \deg(f(t)) > 0 \lor f(t)^{q} = f(t)\right).$$
Note that ``$\deg(f) > 0$'' is universally defined by 9 quantifiers (\cite[Theorem 3.1]{eisenint}) and thus $\fs$ is universally defined in $\fr$ by $\max\{89, 9\} = 89$ quantifiers, as required. \bs


\begin{cor}\label{adjust}
There is a universal $\emptyset$-definition of $\fs$ in $\fr$ given by 90 quantifiers.
\end{cor}

\proof
In the universal definition presented in \emph{Corollary \ref{thecor}}, there are three places parameters are in use: the nonsquare $z \in \mathbb{F}_q$, in $\widetilde{R_{a,b}}$, and in all statements about degree. Examination of \cite[Lemma 3.19]{eisen} reveals $\widetilde{R_{a,b}}$ is defined without use of parameters other than $a$ and $b$, which we already quantify over. To use Eisentr\"{a}ger's formula for degree \cite[Theorem 3.1]{eisenint} without parameters we can define elements of $\mathbb{F}_q \defeq \mathbb{F}_p (\alpha)$ up to conjugates by giving the minimal polynomial for $\alpha$ over $\mathbb{F}_p$. The parameters in Eisentr\"{a}ger's formula are now definable in $\mathcal{L}_{\mbox{\tiny rings}, t} \cup \{\alpha\}$, at the cost of an additional quantifier for $\alpha$. Finally, in this language any nonsquare $z$ is definable, and nonsquares of $\mathbb{F}_q$ remain nonsquare up to conjugates.\bs

\bigskip
\section*{Acknowledgements}
This paper arose from the author's master's thesis, for which Damian R\"{o}ssler was a wonderful supervisor - thank you for that. The author would also like to thank Jochen Koenigsmann for his various insights and assistance along the way, and for his suggestions on \emph{Corollary \ref{adjust}} in particular. Many thanks to Nicolas Daans for sharing his ideas regarding universal definitions of global fields, and for sharing his thoughts on the author's thesis, too. Finally the author extends his thanks to the anonymous referees at Manuscripta Mathematica for their suggested improvements.
\bigskip
\bigskip
\bibliography{refs}   

\begin{thebibliography}{10}

\bibitem{daans}
{\sc Daans, N.}
\newblock {Diophantine definability in number fields and their rings of
  integers}.
\newblock Master's thesis, {Universiteit Antwerpen}, 2018.

\bibitem{daansnew}
{\sc Daans, N.}
\newblock {Universally defining finitely generated subrings of global fields}.
\newblock {\em
  \ttfamily\emph{\href{https://arxiv.org/abs/1812.04372}{arxiv.org/abs/1812.04372}}\/}
  (2018).

\bibitem{dpr}
{\sc Davis, M., Putnam, H., and Robinson, J.}
\newblock {The Decision Problem for Exponential Diophantine Equations}.
\newblock {\em Ann. of Math. (2) 74}, 3 (1961), 425--436.

\bibitem{denef}
{\sc Denef, J.}
\newblock {The Diophantine Problem for polynomial rings of positive
  characteristic}.
\newblock In {\em Studies in Logic and the Foundations of Mathematics},
  M.~Boffa, D.~Dalen, and K.~Mcaloon, Eds., vol.~97. 1979, pp.~131--145.

\bibitem{eisenint}
{\sc Eisentr{\"a}ger, K.}
\newblock {Integrality at a prime for global fields and the perfect closure of
  global fields of characteristic p$>$2}.
\newblock {\em J. Number Theory 114}, 1 (2005), 170--181.

\bibitem{eisen}
{\sc Eisentr{\"a}ger, K., and Morrison, T.}
\newblock {Universally and existentially definable subsets of global fields}.
\newblock {\em Math.\ Res.\ Lett. 25}, 4 (2018), 1173--1204.

\bibitem{hilbert}
{\sc Hilbert, D.}
\newblock {Mathematische Probleme}.
\newblock {\em Nachr.\ K\"{o}nigl.\ Gesell.\ Wiss.\ G\"{o}ttingen,
  Mathematisch-Physikalische Klasse\/} (1900), 253--297.

\bibitem{defzinq}
{\sc Koenigsmann, J.}
\newblock Defining $\mathbb{Z}$ in $\mathbb{Q}$.
\newblock {\em Ann. of Math. 183}, 1 (2016), 73--93.

\bibitem{matiya2}
{\sc Matiyasevich, Y.}
\newblock {Enumerable sets are Diophantine}.
\newblock {\em Soviet Math.\ Dokl. 11\/} (1970), 354--358.

\bibitem{park}
{\sc Park, J.}
\newblock {A universal first order formula defining the ring of integers in a
  number field}.
\newblock {\em Math.\ Res.\ Lett. 20}, 5 (2013), 961--980.

\bibitem{pheidas}
{\sc Pheidas, T.}
\newblock {Hilbert's Tenth Problem for fields of rational functions over finite
  fields}.
\newblock {\em Invent. Math. 103}, 1 (1991), 1--8.

\bibitem{rosen}
{\sc Rosen, M.}
\newblock {\em {Number Theory in Function Fields}}.
\newblock {Springer-Verlag}, 2002.
\newblock {Graduate Texts in Mathematics 210}.

\bibitem{serre79}
{\sc Serre, J.-P.}
\newblock {\em {Local Fields}}.
\newblock Springer-Verlag, 1979.
\newblock {Graduate Texts in Mathematics 67}.

\bibitem{videla}
{\sc Videla, C.}
\newblock {Hilbert's Tenth Problem for rational function fields in
  characteristic 2}.
\newblock {\em Proc. Amer. Math. Soc. 120}, 1 (1994), 249--253.

\end{thebibliography}
\bibliographystyle{acm} 


\end{document}